\documentclass[ijoc,nonblindrev]{informs3} 

\OneAndAHalfSpacedXII 

\usepackage{natbib}
 \bibpunct[, ]{(}{)}{,}{a}{}{,}%
 \def\bibfont{\small}%

 \usepackage{siunitx}
\TheoremsNumberedThrough     

\EquationsNumberedThrough    

\usepackage{booktabs}
\usepackage{pgfplots,tikz,bm}
\pgfplotsset{compat=1.18} 
\usepackage{caption}
\usepackage{subcaption}
\usepackage{multirow}
\usepackage{xcolor}

\newcommand{\pipes}{\mathcal A_{\mathrm{p}}}
\newcommand{\shortpipes}{\mathcal A_{\mathrm{sp}}}
\newcommand{\res}{\mathcal A_{\mathrm{r}}}
\newcommand{\lres}{\mathcal A_{\mathrm{lr}}}
\newcommand{\acelem}{\mathcal A_{\mathrm{active}}}
\newcommand{\comp}{\mathcal A_{\mathrm{c}}}
\newcommand{\cv}{\mathcal A_{\mathrm{cv}}}
\newcommand{\valves}{\mathcal A_{\mathrm{v}}}
\newcommand{\ac}{\mathrm{ac}}
\newcommand{\byp}{\mathrm{bp}}

\usepackage{threeparttable}

\begin{document}



\RUNTITLE{Optimal Gas Flow Problem for a Non-Ideal Gas}

\TITLE{\Large Relaxations of the Steady Optimal Gas Flow \\ Problem for a Non-Ideal Gas}

\ARTICLEAUTHORS{%
\AUTHOR{Sai Krishna Kanth Hari, Kaarthik Sundar, Shriram Srinivasan, and Russell Bent}
\AFF{Los Alamos National Laboratory, \EMAIL{\texttt{\{hskkanth,kaarthik,shrirams,rbent\}@lanl.gov}}}
} 

\ABSTRACT{%
Natural gas ranks second in consumption among primary energy sources in the United States. The majority of production sites are in remote locations, hence natural gas needs to be transported through a pipeline network equipped with a variety of physical components such as compressors, valves, control valves, etc. Thus, from the point of view of both economics and reliability, it is desirable to achieve optimal transportation of natural gas using these pipeline networks. The physics that governs the flow of natural gas through various components in a pipeline network is governed by nonlinear and non-convex equality and inequality constraints and the most general steady-flow operations problem takes the form of a Mixed Integer Nonlinear Program (MINLP). In this paper, we consider one example of steady-flow operations -- the Optimal Gas Flow (OGF) problem for a natural gas pipeline network that minimizes the production cost subject to the physics of steady-flow of natural gas. The ability to quickly determine global optimal solution and a lower bound to the objective value of the OGF for different demand profiles plays a key role in efficient day-to-day operations. One strategy to accomplish this relies on tight relaxations to the nonlinear constraints of the OGF. Currently, many nonlinear constraints that arise due to modeling the non-ideal equation of state either do not have relaxations or have relaxations that scale poorly for realistic network sizes. In this work, we combine recent advancements in the development of polyhedral relaxations for univariate functions to obtain tight relaxations that can be solved \emph{within a few seconds}  on a standard laptop. We demonstrate the quality of these relaxations through extensive numerical experiments on very large scale test networks available in the literature and we find that the proposed relaxation is able to prove optimality in 92\% instances that were used for the experiments. 
}%



\maketitle

\section{Introduction} \label{sec:introduction}

As of 2021, natural  gas ranks second in consumption among primary energy sources in the United States (see report by \cite{EIA}). Cleaner combustion, lower prices, and technological advances have positioned natural gas as the leading source of electricity generation, with its dominance increasing yearly. This ascendancy has garnered the interest of several industries in understanding the economic benefits of optimal inventory planning and operation of natural gas transportation networks. Moreover, this has led to an emergence of a multitude of optimization problem formulations including, but not limited to, expansion planning of natural gas networks (\cite{conradoexpansion}), optimal transportation of natural gas with minimum energy expenditure or production cost (\cite{haristorage}), joint operation of electric and natural gas networks (\cite{roaldgasgrid}), $N$-$k$ interdiction (\cite{ahumadank}), etc. The ability to determine global optimal solutions to these problems is key in adapting quickly to the fluctuations in the demand profiles and prices, and consequently saving millions of dollars. However, the constraints governing the steady-state flow of natural gas in a network are nonlinear and non-convex, thereby rendering the resultant optimization problems computationally challenging to solve. Furthermore, the presence of certain components such as valves and regulators calls for the usage of binary decision variables for modeling. This makes Mixed Integer Nonlinear Programs (MINLPs) the most general form of optimization problems occurring in the context of natural gas networks. It is known that designing convex or polyhedral relaxations for MINLPs that provide tight bounds to the optimal objective value is a key sub-problem for designing global optimization algorithms (see \cite{Floudas2013}). 

In this paper, we restrict our attention to the steady-state Optimal Gas Flow (OGF) problem  for a large scale pipeline network. The dynamics of natural gas flow through pipelines in a regime without waves and shocks can be adequately described by a system of coupled partial differential equations (PDEs) (see \cite{GyryaZlotnik2017}).   
Under steady-state conditions, the system of coupled PDEs reduce to a non-linear system of algebraic equations that relate pressure and mass flow values throughout the pipeline network.  While gas flows in a real pipeline network always undergo a temporal variation, mid-term and long-term planning questions are in practice often examined using steady-state models mainly because, in mid- or long-term planning, future nomination profiles and their time-dependence may not be known. Instead, fictitious future nominations are considered, where load flows and external conditions are assumed to be constant over a fixed time period or the nominations can be aggregated over a day, for example. Furthermore, many optimization problem formulations for gas pipeline networks that feature slow variations over a time horizon simplify flows to sequences of steady-state problems at discrete time instants within the time horizon (see \cite{gugat2021transient}).

The OGF is defined as a pipeline network that consists of components through which gas flows. It is also defined by a set junctions, of which some are in-take (producers) and some are off-take (consumers) points. These are sometimes referred to as delivery points or receipt points, respectively. The flow of gas across different components in a pipeline network is governed by nonlinear, algebraic, steady-state equations that relate the component end-point pressures and flow across the respective component. The OGF problem minimizes total production cost by finding injections at producers to meet the demand from the consumers, while satisfying both the physics of natural gas flow for each component and engineering pressure limits of the network. In particular, the OGF is a MINLP and the focus of this paper is to develop tight relaxations of the non-linear physics of steady-flow for the natural gas using a non-ideal equation of state and demonstrating its efficacy for the original MINLP of interest.

In recent years, most of the research in this area has  focused on developing convex relaxations for the nonlinear, non-convex, and steady-state physics that govern the ideal flow of natural gas in the network (see \cite{conradoexpansion, Singh2019, Singh2020}). The state-of-the-art in this area is the work by authors in \cite{tasseff2020natural} who develop Mixed-Integer Second Order Conic (MISOC) relaxations for many frequently appearing components like valves, resistors, loss resistors and control valves in a gas pipeline network.  Nonetheless, these efforts suffer from  two major shortcomings, namely modeling and scalability. The existing relaxations use an ideal equation of state to relate the pressure and density of the gas at any point in the network; recent studies by authors of \cite{srinivasan2022numerical} have shown that a non-ideal equation of state is more representative of the true nature of flow of natural gas over long pipelines and that the pressure solutions obtained using an ideal and a non-ideal equation of state can vary significantly. In addition, we also show in this paper that the MISOC relaxation is difficult to scale when used on very large scale test networks from a standard benchmark library.
In this paper, we overcome both of these shortcomings by (i) incorporating a non-ideal equation of state and subsequently developing polyhedral relaxations for the gas flow equations in \cite{srinivasan2022numerical} and (ii) introducing models of decision groups which more accurately restrict operational decisions and, as a side effect, improve the scalability of the relaxations on very large pipeline networks. In the process, we also show how the existing MISOC relaxations in \cite{tasseff2020natural} can be modified to work with physics governed by non-ideal equation of state. We present extensive computational experiments to corroborate the effectiveness of all the relaxations in computing global optimal solutions to the OGF. 

The rest of the paper is organized as follows: Section \ref{sec:physics} presents the physics of steady gas flow in a single pipe and their non-dimensional counterparts, followed by the physics for a pipeline network in Section \ref{sec:modeling}. Section \ref{sec:ogf-minlp} -- \ref{sec:misoc-relax} present the full MINLP formulation of the OGF, its linear relaxation, and its MISOC relaxation, respectively. Finally, in Section \ref{sec:results} we present extensive computational results on benchmark instances followed by concluding remarks in Section \ref{sec:conclusion}.

\section{Steady Gas Flow in a Pipe}  \label{sec:physics} 
The adiabatic flow of compressible natural gas in a pipeline is described by the Euler equations in one dimension (see \cite{thorley87}). Ignoring the inertial terms for long pipelines and a steady flow assumption leads to the following equations (see \cite{srinivasan2022numerical}) that represent the mass and momentum conservation equations 
\begin{flalign}
    \varphi \text{ is a constant } \quad \text{ and } \quad \frac{dp}{dx} = -\frac{\lambda}{2D} \frac{\varphi|\varphi|}{\rho} \label{eq:ss-physics}
\end{flalign}
where $\rho$ is the density of the gas, $p$ is pressure, $\phi = \rho v$ is the mass flux, and $v$ is the velocity of the gas. The additional parameters are the friction factor $\lambda$ and diameter $D$ of the pipe. The conservation equations in Eq. \eqref{eq:ss-physics} are supplemented with the equation of state (EoS) relating the $\rho$ and $p$ as 
\begin{flalign}
    \rho(p) = \frac{p}{Z(p, T) R_g T}, \label{eq:eos}
\end{flalign}
where $R_g$ is the specific gas constant, $T$ is the temperature of the gas, $Z$ is the compressibility factor which may in turn depend on the pressure and temperature for a non-ideal gas. It is customary to provide a formula for the compressibility factor $Z(p, T)$ with parameters that have been fitted to measured data obtained using engineering studies (see \cite{benedict1940empirical,elsharkawy2004efficient}). One such widely-used formula for the compressibility factor is the the California Natural Gas Association (CNGA) EoS (see \cite{menon2005gas}) given by 
\begin{flalign}\label{eq:cnga}
    Z(p, T) = \frac 1{b_1 + b_2 p} 
\end{flalign}
where, $b_1$ and $b_2$ are gas and temperature-dependent constants. The values of $b_1$ and $b_2$ are given by the following expressions:
\begin{align}
  b_1 = 1 + \left( \dfrac{p_{atm}}{6894.75729}\right) 
  \left( \dfrac{a_1 10 ^ {a_2 G}}{(1.8T) ^ {a_3}} \right) ~~(\text{\si{unitless}}),\label{eq:b1}\\
    b_2 = \left(\dfrac{1}{6894.75729}\right)
    \left(\dfrac{a_1 10 ^ {a_2 G}}{(1.8T)^{a_3} } \right) ~~(\text{\si{\per\pascal}}). \label{eq:b2}
\end{align}
Here, $b_1$ and $b_2$ are calculated in terms of other non-dimensional constants $a_1 = 344400$, $a_2 = 1.785$, $a_3 = 3.825$, specific gravity of natural gas $G = 288.706$ and atmospheric pressure $p_{atm} = 101350\; \si{\pascal}$. In this context, we remark that for an ideal EoS, the value of $b_1$ and $b_2$ are set to $1$ and $0$, respectively (see \cite{srinivasan2022numerical}). In this paper, we assume isothermal conditions, i.e., $T$ is constant. For the CNGA EoS, the dependence between $\rho$ and $p$ is precisely defined  by combining Eq. \eqref{eq:eos} and \eqref{eq:cnga} as 
\begin{flalign}
    \rho = \frac{b_1 p + b_2 p^2}{R_gT}, \label{eq:rho(p)}
\end{flalign}
where $a > 0$ the fixed quantity $a = \sqrt{R_gT}$. For an ideal EoS the values of $b_1$ and $b_2$ are $1$ and $0$, respectively, and in that case, Eq. \eqref{eq:rho(p)} would correspond to $p = a^2 \rho$.  Rewriting Eq. \eqref{eq:rho(p)} for isothermal conditions using $a$, we obtain\begin{flalign}
    \rho = \frac{b_1 p + b_2 p^2}{a^2}. \label{eq:rho-a}
\end{flalign}
Also, we let $f = A\varphi$ denote the constant mass flow in a pipe with cross-sectional area $A$. Now, integration of Eq. \eqref{eq:ss-physics} along the length of the pipe with end point $1$ and $2$, and with mass flow directed from $1$ to $2$ yields
\begin{align}
    f \text{ is a constant } \quad \text{and } \quad \frac{b_1}2(p_2^2 - p_1^2) + \frac{b_2}3(p_2^3 - p_1^3) 
    = 
    - \dfrac{\lambda La^2}{2DA^2} f |f|\label{eq:ss-final}
\end{align}
where $L$ is the length of the pipe. We note that when the flow is directed from $2$ to $1$, then $f$ is negative. If we define $\pi(p) \triangleq \frac{b_1}2 \cdot p^2 + \frac{b_2}3 \cdot p^3$ to denote the potential at any point along the pipe, then Eq. \eqref{eq:ss-final} can be rewritten as 
\begin{align}
    f \text{ is a constant } \quad \text{and } \quad \pi(p_2) - \pi(p_1)
    = 
    - \dfrac{\lambda La^2}{2DA^2} f |f|\label{eq:ss-final-pi}
\end{align}

\subsection{Non-dimensionalization} 
Previous studies in \cite{srinivasan2022numerical} have shown that non-dimensionalization of the governing equations presented in the previous section plays an important role in convergence of any numerical technique to solve the gas flow equations. To that end, we present the non-dimensionalized equivalent of all the equations presented in the previous section. We refer interested readers to \cite{srinivasan2022numerical} for a detailed discussion on non-dimensionalization and its effect on numerical techniques to solve the steady gas flow equations for a non-ideal gas in a large pipeline network. Here, we let $l_0$, $p_0$, $\rho_0$, $v_0$, $A_0$, $\varphi_0 = v_0 \rho_0$, $f_0 = \varphi_0 A_0$, and $\pi_0 = \rho_0 p_0 a^2$ denote the nominal values for length, pressure, density, velocity, area, mass flux, mass flow, and potential, respectively. Also, we define $\mathcal M = v_0/a$ as the Mach number of the nominal flow velocity and $\mathcal E = p_0 / \left(\rho_0 a^2\right)$ denotes a constant analogous to the Euler number. Using these nominal values and constants, the dimensionless variables are given by
\begin{flalign}
    \bar p = \frac p {p_0},  \quad 
    \bar L = \frac L {l_0},  \quad 
    \bar \rho = \frac \rho {\rho_0},  \quad 
    \bar A = A,  \quad 
    \bar D = \frac D {l_0},  \quad 
    \bar f = \frac f {f_0},  \quad 
    \bar \pi = \frac \pi {\pi_0},  \quad 
\end{flalign}
Using the above notation, non-dimensional EoS and potential is given by  
\begin{flalign}
\bar{\rho} = \bar b_1\bar p + \bar b_2 \bar p^2 \quad \text{ and } \quad \bar \pi(\bar p) = \frac{\bar b_1}2  \bar p^2 + \frac{\bar b_2}3  \bar p^3 \label{eq:eos-cnga}
\end{flalign}
where $\bar b_1 = \mathcal E b_1$ and $\bar b_2 = \mathcal E p_0 b_2$. Finally, the Eq. \eqref{eq:ss-final-pi}, in dimensionless quantities, can be rewritten as 
\begin{align}
    \bar f \text{ is a constant } \quad \text{and } \quad \bar \pi(\bar p_2) - \bar \pi(\bar p_1)
    = 
    - \dfrac{\mathcal M^2}{\mathcal E} \dfrac{\lambda \bar L}{2\bar D\bar A^2} \bar f |\bar f|\label{eq:ss-final-nd}
\end{align}
In the next section, we model the physics of gas flow through various components in a pipeline network. 


\section{Pipeline Network Modeling} \label{sec:modeling} 
A gas pipeline network can be modeled as a graph $\mathcal G = (\mathcal V, \mathcal A)$ where $\mathcal V$ represents the set of nodes and the set $\mathcal A$ denotes the set of arcs that connect any two nodes in the network. In order to facilitate flow of natural gas through the pipeline network, while satisfying its physical limitations, the flow and the pressure in parts of the network are controlled using a subset of arcs in the network. These type of components are referred to as active arcs. An active arc in a network corresponds to a compressor, control valve or a valve. All other arcs are referred to as passive arcs and they correspond to pipes, resistors, loss resistors or short-pipes. 

In the subsequent sections, we present the models for each type of component in a pipeline network. Though models for a subset of components is already known in the literature (see \cite{tasseff2020natural,koch2015evaluating}), we present them again for the sake of completeness.  Throughout the rest of the article, we utilize over-bar to denote dimensionless quantities.

\subsection{Nodes} \label{subsec:nodes}
A node in a gas pipeline network is a junction where one or more arcs meet. Each node $i \in \mathcal V$ of the gas pipeline network is associated with a pressure and a potential variable denoted by $\bar p_i$ and $\bar{\pi}_i$, respectively. Also associated with $i$ is a minimum and maximum allowable operating pressure $\bar{p}_i^{\min}$ and $\bar{p}_i^{\max}$, respectively. The two variables model all the constraints on the node $i$ as follows:
\begin{subequations}
\begin{gather}
    \bar{\pi}_i = \frac{\bar b_1}2  \bar p_i^2 + \frac{\bar b_2}3  \bar p_i^3 \label{eq:node-physics}\\ 
    \bar{p}_i^{\min} \leqslant \bar p_i \leqslant \bar{p}_i^{\max} \label{eq:node-limits}
\end{gather}
\label{eq:node}
\end{subequations}
Eq. \eqref{eq:node-physics} models the relationship between the pressure and the potential at the node $i$ using the non-ideal EoS and Eq. \eqref{eq:node-limits} models the operating limits of the pipeline network at node $i \in \mathcal V$. When using an ideal EoS, Eq. \eqref{eq:node-physics} reduces to $\bar{\pi}_i =  \bar p_i^2/2$ and this is one of the big differences when moving from an ideal to a non-ideal EoS. Next, we present the model for each type of passive arc in the pipeline network.

\subsection{Passive Arcs} \label{subsec:passive}

The four types of passive arcs in any gas pipeline network are (i) pipes, (ii) short pipes, (iii) resistors, and (iv) loss resistors. Here, pipeline (pipes) are the only physical components of the network. A short pipe is a non-physical component that models connections between two geographically co-located nodes; this component is an artifact of bad data issues that is common in large pipeline networks and theoretically, these components can be pre-processed out of the gas pipeline network without changing any of the results. Nevertheless, we refrain from doing so since existing benchmark libraries (see \cite{gaslib}) contain a large number of these components. Finally, resistors and loss resistors are non-physical elements that are used to model variable or fixed pressure loss that is caused either due to bends in the pipes or due to measuring and filtering devices in the pipeline network. In the subsequent sections, we present the model formulations for each of these components. 

\subsubsection{Pipes} -- \label{subsubsec:pipe} 
We let $\pipes \subset \mathcal A$ denote the set of pipes in the network. Each pipe in $\pipes$ connects to nodes $i$ and $j$ is is denoted by $(i, j)$. Associated with each $(i, j) \in \pipes$ is a variable $\bar f_{ij}$ that denotes the constant mass flow rate of the pipe. As a convention, we assume that $\bar f_{ij}$ is positive when the flow is directed from node $i$ to $j$ and negative, otherwise. Analogous to nodes, we let the minimum and maximum mass flow rates through the pipe be denoted by $\bar{f}_{ij}^{\min}$ and $\bar{f}_{ij}^{\max}$, respectively. Also, for notation simplicity we define the resistance of any pipe $(i, j)$ from Eq. \eqref{eq:ss-final-nd} according to the following equation: 
\begin{flalign}
    \beta_{ij} \triangleq \dfrac{\mathcal M^2}{\mathcal E} \dfrac{\lambda \bar L_{ij}}{2\bar D_{ij}\bar A_{ij}^2} \label{eq:pipe-resistance}
\end{flalign}
where, $\bar L_{ij}$, $\bar D_{ij}$, and $\bar A_{ij}$ are the non-dimensional length, diameter and cross-sectional area of the pipe $(i, j)$. 
Then model that governs the flow of gas along the pipe $(i, j) \in \pipes$ is then given by Eq. \eqref{eq:ss-final-nd} and when rewritten using the variables introduced for pipes and nodes is as follows:
\begin{subequations}
\begin{gather}
    \bar \pi_j - \bar \pi_i =  - \beta_{ij} \bar f_{ij} |\bar f_{ij}| \label{eq:pipe-physics} \\ 
    \bar{f}_{ij}^{\min} \leqslant \bar f_{ij} \leqslant \bar{f}_{ij}^{\max} \label{eq:pipe-limits}
\end{gather}
\label{eq:pipe}
\end{subequations}
In Eq. \eqref{eq:pipe-limits}, the value of $\bar{f}_{ij}^{\min}$ may be negative indicating that the flow in the pipe $(i, j)$ can be directed from $j$ to $i$. Eq. \eqref{eq:pipe-physics} is a restatement of Eq. \eqref{eq:ss-final-nd} using notations introduced for a node and a pipe. 
 
\subsubsection{Short Pipes} -- \label{subsec:short-pipe} We let $\shortpipes \subset \mathcal A$ denote the set of short-pipes in the network. As with any arc, any short pipe in the set $\shortpipes$ that connects arbitrary nodes $i, j \in \mathcal N$ is denoted by $(i, j)$. Associated with each short pipe $(i, j) \in \shortpipes$ is a variable $\bar f_{ij}$ that denotes the constant mass flow rate in $(i, j)$. $\bar f_{ij}$ is positive (negative) when the flow is directed from $i$ to $j$ ($j$ to $i$, respectively). The minimum and maximum operating limits for $\bar f_{ij}$ is given by $\bar f_{ij}^{\min}$ and $\bar f_{ij}^{\max}$, respectively. Using these notations, the model that governs the flow of gas along any short pipe $(i, j) \in \shortpipes$ is then given by 
\begin{subequations}
    \begin{gather}
        \bar \pi_j = \bar \pi_i \label{eq:short-pipe-physics} \\ 
        \bar{f}_{ij}^{\min} \leqslant \bar f_{ij} \leqslant \bar{f}_{ij}^{\max} \label{eq:short-pipe-limits}
    \end{gather}
    \label{eq:short-pipe}
\end{subequations}

\subsubsection{Resistors and Loss Resistors} -- \label{subsubsec:resistor} Within a gas pipeline network, many network components like compressors and certain properties of the gas induce a pressure drop along the direction of flow. Some of the common causes of such pressure losses are turbulence due to bends in pipes, filtering and measuring devices, complex piping inside compressor and pressure regulator stations (see \cite{koch2015evaluating}). Most of these effects are highly non-linear and accurate models are absent for most of them. Hence, resistors and loss resistors are used as a surrogate modeling tool for representing these forms of pressure loss in the network. In particular, resistors are usually modelled using one of the two following types: resistors with constant pressure drop, referred to as loss resistors and resistors with a flow-dependent pressure drop, referred to as simply resistors. Both of these type of resistors are modelled as arcs, and we let $\lres \subset \mathcal A$  and $\res \subset \mathcal A$, denote the set of loss resistors and resistors in the pipeline network. 

We first present the constraints that model a loss resistor. To that end, for every $(i, j) \in \lres$ we let $\bar f_{ij}$ denote the constant mass flow rate through the loss resistor and $\bar f_{ij}^{\min}$ and $\bar f_{ij}^{\max}$ denote its operating range. Similar to pipes and short pipes, $f_{ij}$ is allowed to be negative and when it is negative, the flow is directed from node $j$ to node $i$. Also, for each $(i, j) \in \lres$, we let $\Delta \bar p_{ij}$ denote the constant pressure loss that is incurred by the loss resistor along the direction of flow. Then, the constraints that model a loss resistor are given by 
\begin{subequations}
    \begin{gather}
       \bar p_i - \bar p_j = \Delta \bar p_{ij} \cdot \operatorname{sign}\left(\bar f_{ij} \right) \label{eq:lres-physics} \\ 
       \bar{f}_{ij}^{\min} \leqslant \bar f_{ij} \leqslant \bar{f}_{ij}^{\max} \label{eq:lres-limits}
    \end{gather}
    \label{eq:loss-resistor}
\end{subequations}
Eq. \eqref{eq:lres-physics} induces a pressure drop of $\Delta \bar p_{ij}$ along the direction of flow in the loss resistor and when the flow is zero, it enforces the pressure at the nodes $i$ and $j$ to be equal to one another. 

As for resistors, they induce a mass flow rate dependent, pressure drop along the direction of flow through the component. Before, we present the model for a resistor, we introduce relevant notation. We let $\bar f_{ij}$, $\bar f_{ij}^{\min}$, and $\bar f_{ij}^{\max}$ denote the constant mass flow rate through the resistor $(i, j)$, and its minimum and maximum operating limits. The standard pressure-drop model used for a resistor in the literature (see \cite{finnemore2002fluid,lurie2009modeling,koch2015evaluating}) is the Darcy-Weisbach formula with parameters $\zeta_{ij}$ (unit-less drag factor) and $\bar A_{ij}$ (area), given by 
\begin{gather}
    \left( \frac{\zeta_{ij}}{2\bar A_{ij}^2} \frac{\mathcal M^2}{\mathcal E} \right)   \bar f_{ij}^2 = \begin{cases}
    \bar \rho_i (\bar p_i - \bar p_j) & \text{ if flow is from $i \rightarrow j$, i.e., $\bar f_{ij} > 0$} \\ 
    \bar \rho_j (\bar p_j - \bar p_i) & \text{ if flow is from $j \rightarrow i$, i.e., $\bar f_{ij} < 0$} 
\end{cases}
\label{eq:res-standard}
\end{gather}
where, $\rho_i$ and $\rho_j$ are the densities at nodes $i$ and $j$, respectively. 
In this paper, we approximate the standard model in Eq. \eqref{eq:res-standard} further using
\begin{gather}
    \bar \rho_i (\bar p_i - \bar p_j) \; \approx\;  \bar \rho_i \bar p_i - \bar \rho_j \bar p_j \; \approx \; \bar \pi_i - \bar \pi_j \quad \text{ and } \quad \bar \rho_j (\bar p_j - \bar p_i) \; \approx\;  \bar \rho_j \bar p_j - \bar \rho_i \bar p_i \; \approx \; \bar \pi_j - \bar \pi_i \label{eq:approx}
\end{gather}
to obtain simplified model for a resistor as
\begin{subequations}
    \begin{gather}
        \bar \pi_i - \bar \pi_j = \left( \frac{\zeta_{ij}}{2\bar A_{ij}^2} \frac{\mathcal M^2}{\mathcal E} \right)  |\bar f_{ij}| \bar f_{ij} \label{eq:res-physics} \\
        \bar{f}_{ij}^{\min} \leqslant \bar f_{ij} \leqslant \bar{f}_{ij}^{\max} \label{eq:res-limits}
    \end{gather}
    \label{eq:res}
\end{subequations}
\begin{figure}
    \centering
    \includegraphics{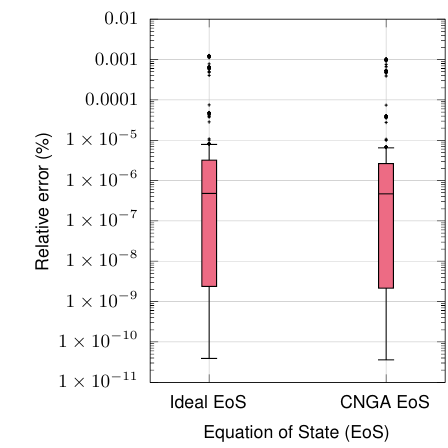}
    \caption{Relative error in outlet pressures obtained using the models in Eq. \eqref{eq:res-standard} and \eqref{eq:res-physics} with the ideal and the CNGA EoS for different values of inlet pressure and flow through the resistors. The ideal EoS is obtained by setting $b_1 = 1$ and $b_2 = 0$ in Eq. \eqref{eq:eos-cnga} (see \cite{srinivasan2022numerical}).}
    \label{fig:res-error}
\end{figure}
The validity of the approximation in Eq. \eqref{eq:res-physics} is verified using the plot shown in Fig. \ref{fig:res-error}; the plot considers $500$ different values of $\bar f_{ij}$ and $\bar p_i$ within its operating limits, for fixed values of $\zeta_{ij}$ and $\bar A_{ij}$ obtained from the GasLib benchmark library (\cite{gaslib}), and compares the relative error between the the values of $\bar p_j$ obtained using models in Eq. \eqref{eq:res-standard} and \eqref{eq:res-physics}, respectively. One added advantage of the new model in Eq. \eqref{eq:res-physics} is that is analogous to the model for a pipe in Eq. \eqref{eq:pipe-physics} with a different value of $\beta_{ij}$. Another appealing reason to use this approximation for resistors is that we can avoid having both density and pressure variables in the formulation of the OGF and the introduction of additional nonlinear constraints to represent their relationship with each other. In the subsequent section, we shall present the models for all the active arcs or the components in the pipeline network.

\subsection{Active Arcs} \label{subsec:active}
The active arcs are controllable arcs of the pipeline network i.e., network operators can either restrict the flow through these components or increase/decrease the pressure across these components as gas flows through the network. We let $\acelem \subset \mathcal A$ denote the subset of all active arcs in the network.  The active arcs in any pipeline network correspond to compressors, control valves, or valves. In the subsequent sections, we present the models for each of these components.  

\subsubsection{Compressors} -- \label{subsubsec:compressor} Compressors are one of the most important components in a gas pipeline network. They increase the pressure of the incoming gas to a higher value thereby enabling transport of gas over long distances by overcoming the friction losses incurred by the passive arcs. Though the internal operation of a compressor is very complex, simple models of its modes of operation is often sufficient for the purpose of solving operational optimization problems (see \cite{tasseff2020natural}). To that end, we let $\comp \subset \acelem$ denote the set of compressors in the pipeline network. Each compressor in $\comp$ connects two nodes $i$ and $j$ and is alternatively denoted by $(i, j)$. Also associated with each compressor $(i, j)$ is a constant mass flow rate $\bar f_{ij}$ which is constrained to be within its operating limits $\bar f_{ij}^{\min}$ and $\bar f_{ij}^{\max}$. Similar to the passive arcs, $f_{ij}$ is positive (negative) when the flow is directed from $i$ to $j$ ($j$ to $i$, respectively). As a convention, for $(i, j) \in \comp$, $i$ is referred to as the inlet node and $j$, the outlet node. Each compressor is also associated with a minimum and maximum compression ratio denoted by $\alpha_{ij}^{\min}$ and $\alpha_{ij}^{\max}$ which limits the amount of compression that can be applied the compressor; typically $\alpha_{ij}^{\min} = 1$. Finally, any compressor can have three modes of operation, namely: \textit{closed}, \textit{bypass}, or \textit{active}. If a compressor is closed, then no gas flows through the compressor and the variables at the inlet and outlet nodes are decoupled. If a compressor is in bypass mode, then the flow direction is dictated by the minimum and maximum values of mass flow that can flow through the compressor but pressure in the inlet node and the outlet node are equal. If a compressor $(i, j)$ is active, then it compressors in the forward direction, flow can be only be directed along the direction of compression and the compression provided by the compressor is within its limits, i.e., $p_j \geqslant p_i$. The conditions for the three modes of operation for any compressor $(i, j)$ is shown in Eq. \eqref{eq:comp-modes}.
\begin{subequations}
    \begin{align}
        \text{closed} \quad & \Rightarrow \quad \text{ $\bar p_i$ and $\bar p_j$ are decoupled} && \text{ and } &&&& \text{$\bar f_{ij} = 0$ } \label{eq:comp-closed} \\
        \text{bypass} \quad & \Rightarrow \quad \text{ $\bar p_i = \bar p_j$} && \text{ and } &&&& \text{$\bar{f}_{ij}^{\min} \leqslant \bar f_{ij} \leqslant \bar{f}_{ij}^{\max}$ } \label{eq:comp-bypass} \\
        \text{active} \quad & \Rightarrow \quad \text{ $\alpha_{ij}^{\min} p_i \leqslant \bar p_j \leqslant \alpha_{ij}^{\max} p_i$} && \text{ and } &&&& \text{$0 \leqslant \bar f_{ij} \leqslant \bar{f}_{ij}^{\max}$ } \label{eq:comp-active}
    \end{align}
    \label{eq:comp-modes}
\end{subequations}

To model the three operating modes in Eq. \eqref{eq:comp-modes}, we introduce three binary variables for each compressor $(i, j) \in \comp$: (i) $x_{ij}$ that takes a value $0$ when the compressor is closed and $1$ when it operates in either the bypass or active mode, (ii) $x_{ij}^{\ac}$ that takes a value of $1$ when compressor is active and $0$ when not, and finally (iii) $x_{ij}^{\byp}$ that takes a value of $1$ when the compressor operates in the bypass mode and $0$, otherwise. Using these variables, the constraints that model all the modes of operation of a compressor $(i, j) \in \comp$ is given by
\begin{subequations}
\begin{gather}
    x_{ij} = x_{ij}^{\ac} + x_{ij}^{\byp} \label{eq:c_status} \\ 
    x_{ij}^{\byp} \cdot \bar f_{ij}^{\min} \leqslant \bar f_{ij} \leqslant x_{ij} \cdot \bar f_{ij}^{\max} \label{eq:c_flow_bounds} \\
    \bar p_j \geqslant \alpha_{ij}^{\min} \bar p_i - (2 - x_{ij}^{\ac} - x_{ij}) \cdot  \left(\alpha_{ij}^{\min} \bar p_i^{\max}  - \bar p_j^{\min}\right) \label{eq:c_p_lb} \\ 
    \bar p_j \leqslant \alpha_{ij}^{\max} \bar p_i + (2 - x_{ij}^{\ac} - x_{ij}) \cdot \left(\bar p_j^{\max} -\alpha_{ij}^{\max} \bar p_i^{\min} \right) \label{eq:c_p_ub} \\ 
    \bar p_i - \bar p_j \geqslant (1 - x_{ij}^{\byp}) \cdot ( \bar p_i^{\min} - \bar p_j^{\max}) \label{eq:c_p_lb_bp} \\ 
    \bar p_i - \bar p_j \leqslant (1 - x_{ij}^{\byp}) \cdot ( \bar p_i^{\max} - \bar p_j^{\min}) \label{eq:c_p_ub_bp} 
\end{gather}
\label{eq:compressor}
\end{subequations}
where Eq. \eqref{eq:c_status} forces the compressor to be in one of the three operating modes, Eq. \eqref{eq:c_flow_bounds} enforces the flow bounds on the compressor for each operating mode, Eq. \eqref{eq:c_p_lb} and \eqref{eq:c_p_ub} enforces the nodal pressure bounds when the compressor is in the active or closed operating mode according to Eq. \eqref{eq:comp-active} and \eqref{eq:comp-closed}, respectively, and finally, Eq. \eqref{eq:c_p_lb} and \eqref{eq:c_p_ub} ensures $\bar p_i = \bar p_j$ when the compressor is the bypass mode. 

\subsubsection{Valves} -- \label{subsubsec:valve} Valves are active arcs in pipeline networks that can be either closed or open. They serve to either route the flow of gas through parts of the network or block the flow completely to parts of the network for maintenance. We let $\valves \subset \acelem$ denote the set of valves in the network and any valve in $\valves$ is denoted by $(i, j)$ where $i$ and $j$ denote the nodes that the valve connects. Associated with each valve $(i, j) \in \valves$ are two variables $\bar f_{ij}$ that is constrained to within its limits $\bar f_{ij}^{\min}$ and $\bar f_{ij}^{\max}$ and $x_{ij} \in \{0, 1\}$ that takes a value $0$ when the valve is closed and $1$, when it is open. Additionally, associated with each valve $(i, j)$ is a parameter $\Delta {\bar p}_{ij}$ that denotes the maximum difference in pressure between $i$ and $j$ when the valve is closed. Given these notations, the constraints that model both the operating modes of a valve are given by 
\begin{subequations}
    \begin{gather}
        x_{ij} \cdot \bar f_{ij}^{\min} \leqslant \bar f_{ij} \leqslant x_{ij} \cdot \bar f_{ij}^{\max} \label{eq:valve_flow_bounds} \\ 
        -\Delta \bar p_{ij} \leqslant \bar p_i - \bar p_j \leqslant \Delta \bar p_{ij} \label{eq:valve_pd} \\ 
        \bar p_i - \bar p_j \geqslant (1 - x_{ij}) \cdot ( \bar p_i^{\min} - \bar p_j^{\max}) \label{eq:valve_open_lb} \\ 
        \bar p_i - \bar p_j \leqslant (1 - x_{ij}) \cdot ( \bar p_i^{\max} - \bar p_j^{\min}) \label{eq:valve_open_ub} 
    \end{gather}
    \label{eq:valve}
\end{subequations}
where, Eq. \eqref{eq:valve_flow_bounds} constrains the flow through the valve to be $0$ when it is closed and to lie between its limits when it is open. Eq. \eqref{eq:valve_open_lb} and \eqref{eq:valve_open_ub} forces the end-point pressures $\bar p_i$ and $\bar p_j$ to be equal to each other when the valve is open and finally, Eq. \eqref{eq:valve_pd} ensures the pressure differential between the ends of the valve is less than the parameter $\Delta \bar p_{ij}$. 

\subsubsection{Control Valves} -- \label{subsubsec:cv} Control valves are the other type of active arc in gas pipeline networks. Contrary to compressors, they reduce the pressure in the control valve's inlet to a lower value and they are usually located at the interface between transmission system and local distribution pipeline networks. Transmission pipeline networks are usually operated at high pressure and distribution network usually operate at a lower pressure and have smaller diameter pipes; the control valves interconnect these two parts of the pipeline network by lowering the pressure to the levels of the distribution network. We let $\cv \subset \acelem$ denote the set of control valves in the systems. Any control valve that connects node $i$ and $j$ is equivalently denoted by $(i, j) \in \cv$. Similar to compressors, they can be operated in three modes: closed, bypass or active. Hence the definitions of the binary variables to model these operating modes and the mass flow variable extend to the control valve as well. Unlike a compressor, any control valve $(i, j)$ has a maximum and minimum pressure differential values $\Delta\bar p^{\max}_{ij}$ and $\Delta\bar p^{\min}_{ij}$ to control the pressure reduction provided and when the control valve is active, the pressure reduction satisfies the constraint $\Delta\bar p^{\min}_{ij} \leqslant \bar p_j - \bar p_i \leqslant \Delta\bar p^{\max}_{ij}$. Using the above notations, the constraints that model the operating modes of any control valve $(i, j)$ are given by 
\begin{subequations}
\begin{gather}
    x_{ij} = x_{ij}^{\ac} + x_{ij}^{\byp} \label{eq:cv_status} \\ 
    x_{ij}^{\byp} \cdot \bar f_{ij}^{\min} \leqslant \bar f_{ij} \leqslant x_{ij} \cdot \bar f_{ij}^{\max} \label{eq:cv_flow_bounds} \\
    \bar p_i - \bar p_j \geqslant (\bar p_i^{\min}  - \bar p_j^{\max}) +  x_{ij}^{\ac} \cdot  \left(\Delta\bar p^{\min}_{ij} - \bar p_i^{\min}  + \bar p_j^{\max}\right) \label{eq:cv_ac_lb} \\ 
    \bar p_i - \bar p_j \leqslant (\bar p_i^{\max}  - \bar p_j^{\min}) - x_{ij}^{\ac} \cdot \left(\bar p_i^{\max} - \bar p_j^{\min} -\Delta\bar p^{\max}_{ij} \right) \label{eq:cv_ac_ub} \\ 
    \bar p_i - \bar p_j \geqslant (1 - x_{ij}^{\byp}) \cdot ( \bar p_i^{\min} - \bar p_j^{\max}) \label{eq:cv_bp_lb} \\ 
    \bar p_i - \bar p_j \leqslant (1 - x_{ij}^{\byp}) \cdot ( \bar p_i^{\max} - \bar p_j^{\min}) \label{eq:cv_bp_ub} 
\end{gather}
\label{eq:control-valve}
\end{subequations}
The description of each constraint in Eq. \eqref{eq:control-valve} is analogous to the constraints that correspond to a compressor in Eq. \eqref{eq:compressor}.

\subsubsection{Sub-network Operation Modes} \label{subsubsec:som} In a large pipeline network, a subset of active arcs, i.e., compressors, valves and control valves can be connected together in more than one way to route gas through the pipeline in different ways. In such cases, modeling the sub-network that consists of these subset of components to reflect all possible ways in which gas can be routed in this sub-network becomes important--it restricts configurations to the small number that are allowable and explicitly introduces cuts on the binary variables that reduce the combinatorial search space. To define sub-network operation modes, we first define a decision group which is a subset of active arcs $\acelem$. Each decision group is associated with one or more operation modes. Each operation mode specifies (a) the on-off status, i.e., \textit{closed} or \textit{open} status, of the component in that group, an optional (b) the direction of flow in each component and (c) optionally for each compressor and control valve in that group, its mode of operation, i.e., \textit{active} or \textit{bypass}. We remark that for each decision group, only one operation mode has to be chosen and all the components in that group have to be operated according to the specifications of that operation mode. 

To model these sub-network operation modes, we use the following notations. We denote by $\zeta$ the set of sub-networks with given operating modes, where each sub-network $S \in \zeta$ is a triple $S \triangleq (\mathcal A_S, \mathcal M_S, f_S)$ of decision group elements $\mathcal A_S \subseteq \acelem$, possible operation modes $\mathcal M_S \subseteq \{0, 1\}^{\mathcal A_S}$ and a function $f_S: \mathcal A_S \times \mathcal M_S \rightarrow \{-1, 0, 1\}$ with 
\begin{gather}
    f_S(a, m) = \begin{cases}
        -1 & \text{if gas flows from $j$ to $i$ in operation mode $m$},\\ 
        0 & \text{if gas flow direction is undefined in operation mode $m$},\\ 
        1 & \text{if gas flows from $i$ to $j$ in operation mode $m$},
    \end{cases} \label{eq:flow-dir}
\end{gather}
for $a = (i, j)$ and $m \in \mathcal M_S$. For each mode $m \in \mathcal M_s$, of a sub-network, we let $s_m \in \{0, 1\}$ indicate if the mode is selected.  For each mode $m$, we let $\mathcal A^{\text{open}}(m)$ and $\mathcal A^{\text{closed}}(m)$ denote the subset of components that are set to open and closed, respectively. Similarly, we let $\mathcal  A^{\ac}(m)$ and $\mathcal A^{\byp}(m)$ denote the subset of control valves and compressors that are operating in active and bypass modes, respectively for $m$. Using these notations, the constraints that capture the behaviour each $S \in \zeta$ are as follows:
\begin{subequations}
    \begin{flalign}
        & \sum_{m \in \mathcal M_S} s_m = 1 \quad \label{eq:one-operating-mode} \\ 
        & s_m \leqslant x_{ij} \quad \forall \;(i, j) \in \mathcal A^{\text{open}}(m), \;m \in \mathcal M_S \label{eq:om-open}\\ 
        & s_m \leqslant 1 - x_{ij} \quad \forall \;(i, j) \in \mathcal A^{\text{closed}}(m), \;m \in \mathcal M_S \label{eq:om-closed}\\
        & s_m \leqslant x_{ij}^{\ac} \quad \forall \;(i, j) \in \mathcal A^{\ac}(m), \;m \in \mathcal M_S \label{eq:om-active}\\ 
        & s_m \leqslant x_{ij}^{\byp} \quad \forall \; (i, j) \in \mathcal A^{\byp}(m), \;m \in \mathcal M_S \label{eq:om-bypass}\\
        & \left( 1 - \sum_{(a, m) \in \mathcal A_S \times \mathcal M_S} s_m \cdot f_S(a, m) \right) \cdot \bar f_{a}^{\min} \leqslant \bar f_a \quad \forall \; a = (i, j) \in \mathcal A_S \label{eq:om-flow-min} \\ 
        & \left( 1 + \sum_{(a, m) \in \mathcal A_S \times \mathcal M_S} s_m \cdot f_S(a, m) \right) \cdot \bar f_{a}^{\max}  \geqslant \bar f_{a} \quad \forall \; a = (i, j) \in \mathcal A_S \label{eq:om-flow-max} 
    \end{flalign}
    \label{eq:operating-modes}
\end{subequations}
where, Eq. \eqref{eq:one-operating-mode} ensures only one operation mode is chosen for each decision group, Eq. \eqref{eq:om-open} -- \eqref{eq:om-closed} sets the on-off status for each component in the decision group based on the operation mode specification, Eq. \eqref{eq:om-active} -- \eqref{eq:om-bypass} sets the compressor and control valves in the decision group to its specification in the operation mode, and finally, Eq. \eqref{eq:om-flow-min} -- \eqref{eq:om-flow-max} enforces flow limits based on the specification in Eq. \eqref{eq:flow-dir}. 

\subsection{Injection and Withdrawal Points} \label{subsec:source-sink} 
Aside from arc components, gas enters or exits the network through injection and withdrawal points located at nodes. Each node $i \in \mathcal V$ may contain more than one injection and withdrawal point. To that end, we let $\mathcal I$ and $\mathcal W$ denote the set of injection and withdrawal points throughout the pipeline network. Additionally, for each node $i \in V$ we let $\mathcal I(i)$ and $\mathcal W(i)$ denote the subset of injections and withdrawals located at node $i$. We assume that the amount of gas that is taken out from the pipeline network at each withdrawal point $i \in \mathcal W$ is known a-priori and is given by $\bar d_i$. Similarly, we let $\bar s_i$ and $\bar s_i^{\max}$ denote the amount of gas that enters the pipeline network at injection point and the maximum injection rate of gas, respectively at injection point $i$. Then, the following constraints are satisfied by the injection and withdrawal points
\begin{subequations}
\begin{gather}
    \sum_{j : (j, i) \in \mathcal A} \bar f_{ji} + \sum_{k \in \mathcal I(i)} \bar s_k = \sum_{j : (i, j) \in \mathcal A} \bar f_{ij} +  \sum_{k \in \mathcal W(i)} \bar d_k \quad \forall i \in \mathcal V \label{eq:nodal-balance} \\
    0 \leqslant \bar s_i \leqslant \bar s_i^{\max} \quad \forall i \in \mathcal I \label{eq:injection-limit} 
\end{gather}
\label{eq:node-equations}
\end{subequations}
where, Eq. \eqref{eq:nodal-balance} enforces the nodal balance at each node in the network and Eq. \eqref{eq:injection-limit} enforces the operation limits on the injection points. Using all the models presented thus far, we are now ready to formulate the OGF problem for a non-ideal gas in a pipeline network as a MINLP. 

\section{Optimal Gas Flow Problem} \label{sec:ogf-minlp}
The OGF problem minimizes the generation cost of gas at the injection points in the pipeline network to meet the gas demand at all the withdrawal points while satisfying constraints that model all the components in the network. In its general form, it is a MINLP. The models for nodes, pipes, resistors and loss resistors require nonlinear equations and models for active arcs require binary variables. To model the generation cost of gas at an injection points, $i$, a linear cost coefficient $c_i$ is used to denote the cost in \$ per unit amount of gas flow rate injected at $i$. Then the OGF is given by:
\begin{equation}
    \tag{MINLP}
    \begin{aligned}[c]
        & \text{minimize:}
        & & \textnormal{Generation cost} && \sum_{i \in \mathcal I} c_i \bar s_i \\
        & \text{subject to:}
        & & \textnormal{Nodal constraints} && ~ \eqref{eq:node}, \eqref{eq:nodal-balance} \\
        & & & \textnormal{Passive arc models} && ~ \eqref{eq:pipe} \text{ -- } \eqref{eq:loss-resistor}, \eqref{eq:res} \\
        & & & \textnormal{Active arc models} && ~ \eqref{eq:compressor} \text{ -- } \eqref{eq:control-valve}, \\
        & & & \textnormal{Sub-network operation modes} && ~ \eqref{eq:operating-modes} \\
        & & & \textnormal{Injection limits} && ~ \eqref{eq:injection-limit}
    \end{aligned}
    \label{eq:minlp}%
\end{equation}
In this MINLP, the constraints in Eq. \eqref{eq:node-physics}, \eqref{eq:pipe-physics}, \eqref{eq:lres-physics} and \eqref{eq:res-physics} are the sources of nonlinearities. Before we present linear and second-order conic relaxations for the rest of the nonlinear terms in the Eq. \eqref{eq:node-physics}, \eqref{eq:pipe-physics}, and \eqref{eq:res-physics}, we first present an exact reformulation of the nonlinear constraint in Eq. \eqref{eq:lres-physics} by the introduction of a binary variable to model the direction of flow of gas in a loss resistor. 

\subsection{Reformulation of the Loss Resistor Model} \label{subsec:lr-reformulation} The nonlinearity in Eq. \eqref{eq:lres-physics} is introduced with the `sign` function. As the `sign` function is non-differentiable, it often cannot be handled with off-the-shelf commercial or open-source nonlinear programming solvers. To reformulate the `sign` function, for each loss resistor $(i, j) \in \lres$, a binary variable $x_{ij}$ is introduced that takes a value of $1$ if the gas flow in the loss resistor is directed from node $i$ to node $j$ and $0$, otherwise. Then, Eq. \eqref{eq:loss-resistor} for loss resistor $(i, j) \in \lres$ is equivalently reformulated as:
\begin{gather}
    \bar{p}_i - \bar{p}_j = \Delta \bar p_{ij} \cdot (2 x_{ij} - 1) \quad \text{and} \quad (1-x_{ij}) \cdot \bar f^{\min} \leqslant \bar f_{ij} \leqslant x_{ij} \cdot \bar f^{\max} \label{eq:loss-resistor-eq}
\end{gather}
In Eq. \eqref{eq:loss-resistor-eq}, it is assumed that $\bar f^{\min} <0$, i.e., the flow through a loss resistor can either be directed from $i$ to $j$ or from $j$ to $i$. If this is not the case, then the model for a loss resistor is fully linear and this reformulation is unnecessary\footnote{The absolute value function in Eq. \eqref{eq:pipe-physics} can be reformulated in a similar way \cite{conradoexpansion}}. 

\section{Linear Relaxation} \label{sec:lin-relax}
In this section, we present one of the key contributions of this paper, a linear relaxations for the nonlinear terms of Eq. \eqref{eq:node-physics}, \eqref{eq:pipe-physics}, and \eqref{eq:res-physics}. To present the relaxations in a concise manner, we formulate a linear relaxations for arbitrary univariate functions that are continuous, have a bounded domain, and are differentiable in their domain, i.e., $y = g(x)$ where $g : [a, b] \rightarrow \mathbb R$ and $g \in \mathcal C^1$ where $\mathcal C^1$ denotes the space of continuous and once differentiable functions. For Eq. \eqref{eq:node-physics} the function $g(x) =  (\bar b_1/2)  x^2 + (\bar b_2/3)  x^3$ (with the domain $[\bar p_i^{\min}, \bar p_i^{\max}]$ for each node $i \in \mathcal V$), and for Eq. \eqref{eq:pipe-physics} and \eqref{eq:res-physics}, $g(x) = x|x|$ (with the domain $[\bar f_{ij}^{\min}, \bar f_{ij}^{\max}]$ for each $(i, j) \in \pipes \cup \res$). It is clear that, in both cases, $g$ is univariate and in $\mathcal C^1$. Given this abstraction, we construct a linear relaxation of the constraint $y = g(x)$ with $x \in [a, b]$. Here we leverage recent results in obtaining a sequence of polyhedral relaxations for univariate functions in \cite{sundar2021sequence}. For the sake of completeness, we present relevant definitions that are useful for presenting the linear relaxation and provide a geometric intuition for these relaxations. 

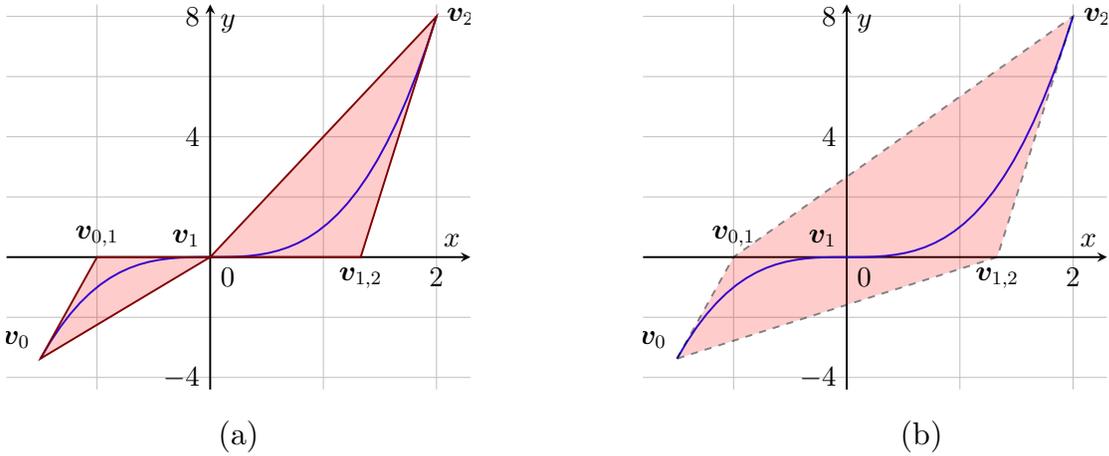
\begin{figure}[htbp]
\centering
\begin{subfigure}[b]{0.45\textwidth}
\centering
\begin{tikzpicture}[scale=0.9]
\begin{axis}[ xlabel={$x$}, ylabel={$y$}
  ,axis lines=middle
  ,xmin=-1.8, xmax=2.3, ymin=-4.4, ymax=8.4
  ,samples=41, grid, thick
  ,domain=-1.5:2
  ,legend pos=outer north east
  ,xtick={}
  ,ytick={}
  ,ticks=none
  ]
  \addplot+[no marks] {x^3};
  \addplot [red!50!black,fill=red,fill opacity=0.2]
  coordinates 
  {(0, 0) (2, 8) (1.33, 0.0)} -- cycle;
  \addplot [red!50!black,fill=red,fill opacity=0.2]
  coordinates 
  {(-1.5, -3.375) (0, 0) (-1, 0.0)} -- cycle;
  \node [above left] at (-1.5, -3.375) {$\bm v_0$};
  \node [above left] at (0, 0) {$\bm v_1$};
  \node [right] at (2, 8) {$\bm v_2$};
  \node [above] at (-1, 0.0) {$\bm v_{0, 1}$};
  \node [below] at (1.33, 0.0) {$\bm v_{1, 2}$};
  \node [below right] at (0.0, 0.0) {$0$};
  \node [below] at (2, 0) {$2$};
  \node [left] at (0, 4) {$4$};
  \node [left] at (0, 8) {$8$};
  \node [left] at (0, -4) {$-4$};
 \end{axis}
\end{tikzpicture}
\caption{}
\label{fig:x3-triangles}
\end{subfigure} 
\hfill 
\begin{subfigure}[b]{0.45\textwidth}
    \begin{tikzpicture}[scale=0.9]
    \begin{axis}[xlabel={$x$}, ylabel={$y$}
  ,axis lines=middle
  ,xmin=-1.8, xmax=2.3, ymin=-4.4, ymax=8.4
  ,samples=41, grid, thick
  ,domain=-1.5:2
  ,legend pos=outer north east
  ,xtick={}
  ,ytick={}
  ,ticks=none
  ]
  \addplot+[no marks] {x^3};
  \addplot [dashed,opacity=0.5,fill=red,fill opacity=0.2]
  coordinates {(-1.5, -3.375) (-1, 0.0) (2, 8) (1.33, 0.0)} -- cycle;
  \node [above left] at (-1.5, -3.375) {$\bm v_0$};
  \node [above left] at (0, 0) {$\bm v_1$};
  \node [right] at (2, 8) {$\bm v_2$};
  \node [above] at (-1, 0.0) {$\bm v_{0, 1}$};
  \node [below] at (1.33, 0.0) {$\bm v_{1, 2}$};
  \node [below right] at (0.0, 0.0) {$0$};
  \node [below] at (2, 0) {$2$};
  \node [left] at (0, 4) {$4$};
  \node [left] at (0, 8) {$8$};
  \node [left] at (0, -4) {$-4$};
 \end{axis}
\end{tikzpicture}
\caption{}
\label{fig:x3-relaxation}
\end{subfigure}
\caption{The construction of the linear relaxations for the function $y = x^3$ with domain $[-1.5, 2]$. Here, $p_g^{\text{base}} = \{-1.5, 0, 2\}$. The coordinates of the vertices are $\bm v_0 = (-1.5, -3.375)$,  $\bm v_1 = (0, 0)$,  $\bm v_2 = (2, 8)$,  $\bm v_{0, 1} = (-1, 0)$,  $\bm v_{1, 2} = (1.33, 0)$. (a) shows the triangles corresponding to each sub-interval of $p_g^{\text{base}}$ and (b) shows the convex hull of the triangles which gives the linear relaxation of $y = x^3$. Observe in (a) that for each sub-interval the curve is completely contained in the triangle for that sub-interval. In (b), $\bm v_{0, 1}$ and $\bm v_{1, 2}$ are the intersection of the tangents for each sub-interval and the shaded region gives the convex hull of the triangles in (a).}
\label{fig:x3-lp}
\end{figure}

Here, we define a partition of the domain of $g$ as the set $p \triangleq \{a = x_0, x_1, \dots, x_{n-1}, x_n = b\}$ with $x_0 < x_1 < \cdots < x_{n-1} < x_n$. Also, we define a base partition of $g$ (denoted by $p_g^{\text{base}}$) as any partition that satisfies the following conditions: (i) all the break-points (the points in the domain where the function changes from being convex to concave or vice-versa) of the function $g$ are contained in $p_g^{\text{base}}$, (ii) for any successive partition points the derivative values are not equal to each other, i.e, $g'(x_i) \neq g'(x_{i+1})$ and finally, (iii) the cardinality of $p_g^{\text{base}}$ is minimum. For any node $i \in \mathcal V$, the base partition of $g(x) =  (\bar b_1/2)  x^2 + (\bar b_2/3)  x^3$ is $p_g^{\text{base}} = \{\bar p_i^{\min}, \bar p_i^{\max}\}$ since $\bar p_i^{\min} > 0$ and the function is convex in its domain and for any $(i, j) \in \pipes \cup \res$, the base partition of $g(x) = x|x|$ is $p_g^{\text{base}} = \{\bar f_{ij}^{\min}, 0, \bar f_{ij}^{\max}\}$ when $\bar f_{ij}^{\min} < 0$ and $\bar f_{ij}^{\max} > 0$  and $p_g^{\text{base}} = \{\bar f_{ij}^{\min}, \bar f_{ij}^{\max}\}$, otherwise. Given, $y = g(x)$ and $p_g^{\text{base}}$, the linear relaxation for an arbitrary $g(x)$ with domain $[a, b]$ is is now constructed using the following procedure:
\begin{enumerate}
    \item Let $y = g(x)$ and $p_g^{\text{base}} = \{a = x_0, x_1, \dots, x_{n-1}, x_n = b\}$. 
    \item For each sub-interval $[x_i, x_{i + 1}]$ with $i \in \{0, \cdots, n-1\}$, construct a triangle that is defined by the two tangents to $y = g(x)$ at $x_i$ and $x_{i+1}$ and the secant line from $g(x_i)$ to $g(x_{i+1})$. We observe that for the sub-interval $[x_i, x_{i + 1}]$, this triangle is a relaxation of the graph of $y = g(x)$ restricted to that sub-interval; this is always guaranteed because within that sub-interval the definition of $p_g^{\text{base}}$ ensures that the restriction of $g$ to $[x_i, x_{i + 1}]$ is either convex or concave. See Fig. \ref{fig:x3-triangles} for a geometric illustration. 
    \item The linear relaxation of $y = g(x)$ in its domain $[a, b]$ is then given by the convex hull of all the triangles, one in each sub-interval. See Fig. \ref{fig:x3-lp} for a geometric illustration. 
\end{enumerate}

The relaxations obtained by this procedure can be tightened by further refining the base partition $p_g^{\text{base}}$ i.e., by including more partition points in $p_g^{\text{base}}$. It is known that for a certain classes of refinement schemes (see \cite{sundar2021sequence}) that successively refine the base partition $p_g^{\text{base}}$, the linear relaxation converges to the convex hull of the $y = g(x)$. Now, it remains to mathematically characterize the convex hull of the triangles obtained using this constructive procedure. We will do so for the base partition $p_g^{\text{base}}$ which generalizes to any refinement of the base partition. Here, $\bm v_0$ and $\bm v_n$ denotes vertices corresponding to the points $(x_0, g(x_0))$ and $(x_n, g(x_n))$. Then, $\bm v_{i, i+1}$ denotes the vertex of the triangle for the sub-interval $[x_i, x_{i+1}]$ that is the intersection of the two tangents at $x_i$ and $x_{i+1}$. The mathematical characterization of the convex hull of the triangles, denoted by $\left\{ y = g(x)\right \}^{\text{LP}}$, is calculated using the observation that the only vertices of a triangle that constitute the extreme points of convex hull of all triangles are in the set $\mathcal K = \{\bm v_0, \bm v_n\} \cup \{\bm v_{0, 1}, \dots, \bm v_{k-1, k}\}$ where $k = n-1$. This is because the other vertices can be expressed as a convex combination of vertices in $\mathcal K$. In fact, when $g$ is convex or concave, $\mathcal K$ is exactly the set of extreme points of convex hull of the triangles. Hence, a linear programming formulation $\left\{ y = g(x) \right \}^{\text{LP}}$ is given by the following convex hull description of vertices in $\mathcal K$:
\begin{flalign}
\left\{ y = g(x) \right \}^{\text{LP}} \triangleq \left\{(x, y, \bm \lambda) \in [a, b] \times \mathbb R \times \bm \Delta_{|\mathcal K|} : \begin{pmatrix} x \\ y \end{pmatrix} = \sum_{i=1}^{|\mathcal V|} \lambda_i \bm w_i \right\} \label{eq:lp}
\end{flalign}
where, $\bm w_1, \bm w_2, \dots, \bm w_{|\mathcal K|}$ are vertices in $\mathcal K$, and $\bm \Delta_{|\mathcal K|}$ is a $|\mathcal K|$-dimensional simplex. Note that $\left\{ y = g(x) \right \}^{\text{LP}}$ Eq. \eqref{eq:lp} is also a function of $p_g^{\text{base}}$ since it is the partition that defines the triangles and in turn the extreme points in $\mathcal K$. This relaxation is applied to all three nonlinear terms Eq. \eqref{eq:node-physics}, \eqref{eq:pipe-physics}, and \eqref{eq:res-physics} to obtain a linear relaxation of the MINLP in the Sec. \ref{sec:ogf-minlp}. In particular, for the Eq. \eqref{eq:node-physics} is relaxed as
\begin{gather}
    \left\{\bar \pi_i = g(\bar p_i)\right\}^{\text{LP}} ~~ \text{where $g: [\bar p_i^{\min}, \bar p_i^{\max}] \rightarrow \mathbb R$, $g(x) = \frac{\bar b_1}2  x^2 + \frac{\bar b_2}3  x^3$} \quad \forall i \in \mathcal V \label{eq:lp-node-physics}
\end{gather}
Similarly, for Eq. \eqref{eq:pipe-physics} and \eqref{eq:res-physics}, the relaxations are 
\begin{subequations}
\begin{gather}
    \bar \pi_j - \bar \pi_i =  - \beta_{ij} \bar F_{ij} \quad \forall (i, j) \in \pipes \label{eq:pipe-physics-lifted} \\ 
     \bar \pi_i - \bar \pi_j = \left( \frac{\zeta_{ij}}{2\bar A_{ij}^2} \frac{\mathcal M^2}{\mathcal E} \right)  \bar F_{ij} \quad \forall (i, j) \in \res \label{eq:res-physics-lifted} \\ 
    \left\{\bar F_{ij} = g(\bar f_{ij})\right\}^{\text{LP}} ~~ \text{where $g: [\bar{f}_{ij}^{\min}, \bar{f}_{ij}^{\max}] \rightarrow \mathbb R$, $g(x) = x|x|$} \quad \forall (i,j) \in \pipes \cup \res \label{eq:lp-piperes-physics}
\end{gather}
\label{eq:lp-pipe-res}
\end{subequations}
In Eq. \eqref{eq:pipe-physics-lifted} and \eqref{eq:res-physics-lifted}, $\bar F_{ij}$ is an auxiliary variable that is used to isolate the nonlinear terms in Eq. \eqref{eq:pipe-physics} and \eqref{eq:res-physics}, respectively. Together, the linear relaxation (LR) for the MINLP is stated as 
\begin{equation}
    \tag{LR}
    \begin{aligned}[c]
        & \text{minimize:}
        & & \textnormal{Generation cost} && \sum_{i \in \mathcal I} c_i \bar s_i \\
        & \text{subject to:}
        & & \textnormal{Nodal constraints} && ~ \eqref{eq:node-limits}, \eqref{eq:nodal-balance}, \eqref{eq:lp-node-physics} \\
        & & & \textnormal{Passive arc models} && ~ \eqref{eq:pipe-limits}, \eqref{eq:short-pipe}, \eqref{eq:res-limits}, \eqref{eq:loss-resistor-eq}, \eqref{eq:lp-pipe-res} \\
        & & & \textnormal{Active arc models} && ~ \eqref{eq:compressor} \text{ -- } \eqref{eq:control-valve} \\
        & & & \textnormal{Sub-network operation modes} && ~ \eqref{eq:operating-modes} \\
        & & & \textnormal{Injection limits} && ~ \eqref{eq:injection-limit}
    \end{aligned}
    \label{eq:milp}%
\end{equation}

In the next section, we present a MISOC relaxation for the nonlinear terms in \eqref{eq:pipe-physics} and \eqref{eq:res-physics}, motivated by existing work by authors in \cite{conradoexpansion,tasseff2020natural,Singh2020} as a state-of-the-art comparison point for our linear relaxations.

\section{Mixed Integer Second Order Cone Relaxation} \label{sec:misoc-relax}
In this section, we present a convex relaxation for the MINLP formulation of the OGF that extends existing MISOC relaxation in the literature (see \cite{conradoexpansion}) to the non-ideal EoS setting. We relax the nonlinear term that arises in the pipe and resistor constraints in \eqref{eq:pipe-physics} and \eqref{eq:res-physics} using the existing MISOC relaxation in the literature. As for the relaxation of the Eq. \eqref{eq:node-physics}, we reuse the linear relaxation developed in the previous section. The combination of both the relaxations applied to the corresponding terms enable direct application of off-the-shelf commercial and open-source mixed-integer convex optimization solvers to compute the optimal solution. Though the MISOC relaxation is well-known in the literature, we present the relaxation in this section to keep the paper self-contained. To start, Eq. \eqref{eq:pipe-physics} and \eqref{eq:res-physics} are reformulated by introducing a binary variable $x_{ij}$ for each $(i, j) \in \pipes \cup \res$ that takes a value $1$ when the flow is directed from $i$ to $j$ ($\bar f_{ij} \geqslant 0$) and $0$, when flow is directed from $j$ to $i$ ($\bar f_{ij} \leqslant 0$). Using this variable, Eq. \eqref{eq:pipe-physics} and \eqref{eq:res-physics} is equivalently reformulated as
\begin{subequations}
\begin{gather}
    \bar \pi_j - \bar \pi_i =  - \beta_{ij}  (2 x_{ij} - 1) \bar f_{ij}^2 \quad \forall (i, j) \in \pipes \label{eq:pipe-physics-eq} \\ 
     \bar \pi_i - \bar \pi_j = \left( \frac{\zeta_{ij}}{2\bar A_{ij}^2} \frac{\mathcal M^2}{\mathcal E} \right)   (2 x_{ij} - 1) \bar f_{ij}^2 \quad \forall (i, j) \in \res \label{eq:res-physics-eq} 
\end{gather}
\label{eq:pipe-res-physics-eq}
\end{subequations}
The nonlinear term in both the above equations is $(2 x_{ij} -1) \bar f_{ij}^2$ which is in-turn can equivalently be written as 
\begin{gather}
    \gamma_{ij} = (2 x_{ij} -1) \hat f_{ij} \;\; \text{ and } \;\; \hat f_{ij} = \bar f_{ij}^2 \quad \forall (i, j) \in \pipes \cup \res \label{eq:lift} 
\end{gather}
where, $\hat f_{ij}$ and $\gamma_{ij}$ are additional auxiliary variables introduced for each $(i, j) \in \pipes \cup \res$. Now, the first equation is \eqref{eq:lift} is relaxed using a standard McCormick relaxation (see \cite{mccormick1976computability}) which, in this case, is exact and the second equation is relaxed as $\hat f_{ij} \geqslant \bar f_{ij}^2$, a second order cone constraint. The McCormick relaxation for $\gamma_{ij} = (2 x_{ij} -1) \hat f_{ij}$ can be constructed using upper and lower bounds for $(2 x_{ij} -1)$ and $\hat f_{ij}$ which is  given by $[-1, +1]$ and $[0, \max \left( |\bar f_{ij}^{\min}|, |\bar f_{ij}^{\max}| \right)^2]$, respectively. Now, the McCormick relaxation for every $(i, j) \in \pipes \cup \res$ is given by
\begin{subequations}
    \begin{gather}
        -\hat f_{ij} \leqslant \gamma_{ij} \leqslant \hat f_{ij} \\ 
        \gamma_{ij} \geqslant \hat f_{ij} + (2 x_{ij} - 1) \cdot \max \left( |\bar f_{ij}^{\min}|, |\bar f_{ij}^{\max}| \right)^2 - \max \left( |\bar f_{ij}^{\min}|, |\bar f_{ij}^{\max}| \right)^2 \\ 
        \gamma_{ij} \leqslant -\hat f_{ij} + (2 x_{ij} - 1) \cdot \max \left( |\bar f_{ij}^{\min}|, |\bar f_{ij}^{\max}| \right)^2 + \max \left( |\bar f_{ij}^{\min}|, |\bar f_{ij}^{\max}| \right)^2
    \end{gather}
    \label{eq:mccormick}
\end{subequations}
Now, the MISOC relaxation for the MINLP in Sec. \ref{sec:ogf-minlp} is obtained by replacing Eq. \eqref{eq:lp-pipe-res}  in LR by the relaxation described in this section. 

\section{Computational Results} \label{sec:results} 
In this section, we present results of extensive computational experiments performed on benchmark gas pipeline network instances. All the presented formulations were implemented in the Julia Programming language (see \cite{bezanson2017julia}) using JuMP (see \cite{DunningHuchetteLubin2017}) as the mathematical programming layer. Furthermore, the Julia package ``PolyhedralRelaxations'' in \url{https://github.com/sujeevraja/PolyhedralRelaxations.jl} was used to formulate the linear relaxations for univariate functions in Sec. \ref{sec:lin-relax}; given any arbitrary univariate $\mathcal C^1$ function and a partition of the domain of the function with all break-points included. All computational experiments were run on an Intel Haswell 2.6 GHz, 62 GB, 20-core machine at Los Alamos National Laboratory. Furthermore, CPLEX, a commercial MILP and MISOC solver, was used to solve all the the relaxations and SCIP (see \cite{achterberg2009scip}) was used to solve the MINLP formulation of the OGF. A computational time limit of 1000 seconds was set for every solve of the MINLP, linear and MISOC relaxations. The code to reproduce all the results presented in the subsequent sections is open-sourced and made available at \url{https://github.com/kaarthiksundar/GasSteadyOpt.jl}. We now present all the data sources and the modifications that were done for the OGF problem considered in this article. 

\subsection{Data Sources} \label{subsec:data-sources}
The computational experiments were performed on gas network instances obtained from GasLib, a library of gas network instances from \cite{gaslib}. In particular, we test the efficacy of our formulations on three networks: GasLib-134, GasLib-582, GasLib-4197 that were generated using real pipeline system data in Europe. The number of network components in each of these pipeline networks is given by Table \ref{tab:gaslib} 
\begin{table}[]
    \centering
        \caption{Data on the size of each network in GasLib}
    \label{tab:gaslib}
    \begin{tabular}{ccccccccc}
        \toprule 
        instance & $|\mathcal V|$ & $|\pipes|$ & $|\comp|$ & $|\cv|$ & $|\res \cup \lres|$ & $|\valves|$ & $|\shortpipes|$ & $|\zeta|$ \\ 
        \midrule 
        GasLib-4197 & 4197 & 3537 & 12 & 120 & 28 & 426 & 343 & 413 \\
        GasLib-582 	& 582 & 278 & 5 & 23 & 8 & 26 & 269 & 2 \\
        GasLib-134 	& 134 & 86 & 1 & 1 & 0 & 0 & 45 & 0 \\
        \bottomrule 
    \end{tabular}
\end{table}
For each of these networks, multiple nomination cases i.e., multiple withdrawal cases are provided. In particular, there are $1234$, $4227$, and $2859$ instances for GasLib-134, GasLib-582, and GasLib-4197, respectively. As a part of these nomination cases, the data set also provides injection values. We set the maximum injection for each injection point to be 5\% higher than  the given maximum. Finally, the  cost in \$ per unit amount of gas flow rate injected at each injection point was assigned a randomly generated value in the range of $[1, 5]$ units. The efficacy of all the formulations were tested on a set of $8320$ instances. 

\subsection{Results for GasLib-134} \label{subsec:gaslib134}
Out of the $1234$ instances for GasLib-134, in $2$ instances both the linear and MISOC relaxations were infeasible implying that the MINLP is also infeasible. For the remaining $1232$ instances, the statistics of the computation time and the relative gaps are given in Table \ref{tab:gaslib-134-stats}. Throughout the rest of the article, the relative gap is measured between the objective value of the relaxation and the objective value of the feasible solution provided by the MINLP solver, relative the relaxation's objective. From the table, it is clear that the relaxations are tight for the feasible MINLP instances in GasLib-134. Nevertheless, it is computationally tractable to directly solve the MINLP for these networks as they can be solved in less than a second. The essential take-away is that for small networks, it is good practice to directly solve the MINLP before relying on relaxations. 

\begin{table}[h!]
    \centering
    \caption{Statistics of computation times and relative gaps for the 1234 feasible GasLib-134 instances.}
    \label{tab:gaslib-134-stats}
    \begin{tabular}{cccccc}
         \toprule 
         \multirow{2}{*}{statistic}& \multicolumn{3}{c}{time (sec.)} & \multicolumn{2}{c}{rel. gap (\%)} \\ 
         \cmidrule(lr){2-4} \cmidrule(lr){5-6}
         & MINLP & MISOC & linear & MISOC & linear \\ 
         \midrule 
         minimum & 0.01 & 0.02 & 0.02 & 0.00 & 0.00 \\
         maximum & 0.12 & 0.55 & 0.06 & 0.00 & 0.00\\ 
         mean & 0.05 & 0.02 & 0.02 & 0.00 & 0.00\\ 
         std. dev. & 0.01 & 0.02 & 0.01 & 0.00 & 0.00\\   
         \bottomrule
    \end{tabular}
\end{table}

\subsection{Results for GasLib-582} \label{subsec:gaslib582}
For the GasLib-582 networks, we start to see the MINLP struggling to solve the OGF, where as the linear and the MISOC relaxations are very effective in computing lower bounds within a few second of computation time. Table \ref{tab:gaslib-582-counts} shows the total number of feasible, infeasible and timed out instances for the MINLP, linear and MISOC relaxations. We observe from the table that there are $7$ instances which are detected to be infeasible by the MINLP, but the linear and MISOC relaxation declare that the relaxations are feasible. This provides an empirical proof that the neither relaxation is tight for the OGF problem. For these $7$ instances SCIP detected infeasibility and both relaxations were solved to optimality with a second of computation time. On the other hand, for the $5$ instances that MINLP times out after the computation time limit of $1000$ seconds, both relaxations converged to their respective optimal solutions in less than a second of computation time and hence, we cannot compute relative gaps. 

\begin{table}[h!]
    \centering
    \caption{Number of instances solved to optimality, infeasible and errored out instances in the GasLib-582 network.}
    \label{tab:gaslib-582-counts}
    \begin{tabular}{ccccc}
        \toprule 
        formulation & \# instances & \# optimal & \# infeasible & \# time limit \\ 
        \midrule 
        MINLP & 4227 & 4215 & 7 & 5 \\ 
        Linear relaxation & 4227 & 4227 & 0 & 0  \\
        MISOC relaxation & 4227 & 4227 & 0 & 0 \\
        \bottomrule
    \end{tabular}
\end{table}

For the $4215$ instances that the MINLP computed the globally optimal solution to the OGF, the statistics of the computation time and the relative gaps are given in Table \ref{tab:gaslib-582-stats}. Overall, the linear and the MISOC relaxations provided a proof of global optimality for $4196/4215$ and $4212/4215$ instances, respectively. Based on the results in Table \ref{tab:gaslib-582-stats}, for medium-sized pipeline networks, we see that the linear relaxation is slightly faster than the MISOC relaxation and the MINLP struggles on a few instances. Hence, the recommendation for the medium-sized instances is to run the MINLP with a small computation time limit to obtain a feasible solution and solve one of the relaxations to obtain an estimate of the relative gap using the cost of the feasible solution and the objective value of the relaxation. 

\begin{table}[h!]
    \centering
    \caption{Statistics of computation times and relative gaps for the 4215 feasible GasLib-582 instances.}
    \label{tab:gaslib-582-stats}
    \begin{tabular}{cccccc}
         \toprule 
         \multirow{2}{*}{statistic}& \multicolumn{3}{c}{time (sec.)} & \multicolumn{2}{c}{rel. gap (\%)} \\ 
         \cmidrule(lr){2-4} \cmidrule(lr){5-6}
         & MINLP & MISOC & linear & MISOC & linear \\ 
         \midrule 
         minimum & 0.26 & 0.12 & 0.15 & 0.00 & 0.00 \\
         maximum & 667.35 & 5.36 & 0.44 & 0.61 & 0.97\\ 
         mean & 8.16 & 0.36 & 0.25 & 0.00 & 0.00\\ 
         std. dev. & 22.27 & 0.15 & 0.06 & 0.00 & 0.00\\   
         \bottomrule
    \end{tabular}
\end{table}

\subsection{Results for GasLib-4197} \label{subsec:gaslib4197}
The GasLib-4197 network is a very large network and we remark that, to the best of our knowledge, \emph{no computational study has been preformed previously in any form of the gas flow problem for networks of this size}. For this network, the total number of instances is $2859$. Table \ref{tab:gaslib-4197-counts} shows the total number of instances that were solved to optimality, number of instances that were reported infeasible, number of instances that timed out within the computation time limit of $1000$ seconds and the number of instances for which the solver errors out. From the table, it is clear that the linear relaxation clearly outperforms the MISOC relaxation since there are $203$ instances for which the MISOC relaxation times out. 
\begin{table}[h!]
    \centering
    \caption{Number of instances solved to optimality, infeasible, errored out and timed-out instances in the GasLib-4197 network.}
    \label{tab:gaslib-4197-counts}
    \begin{tabular}{cccccc}
        \toprule 
        formulation & \# instances & \# optimal & \# infeasible & \# time limit & \# error\\ 
        \midrule 
        MINLP & 2859 & 2225 & 68 & 561 & 5 \\ 
        Linear relaxation & 2859 & 2859 & 0 & 0 & 0 \\
        MISOC relaxation & 2859 & 2656 & 0 & 203 & 0 \\
        \bottomrule
    \end{tabular}
\end{table}

\begin{figure}
    \centering
    \includegraphics{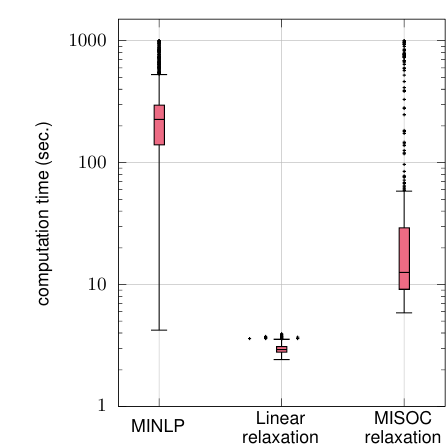}
    \caption{Box plot of computation times for the instances of the GasLib-4197 network were solved by all three formulations.}
    \label{fig:4197-times-success}
\end{figure}

The Fig. \ref{fig:4197-times-success} shows the box plot of the computation times using MINLP, linear relaxation and the MISOC relaxation on the instances for which all three formulations were successful in computing the respective optimal solutions within the computation time limit (a total of $2058$ instances). The plot clearly shows the computational superiority of the linear relaxation compared to the MISOC relaxation; furthermore, the plot also indicates that for many instances of GasLib-4197, \emph{the MISOC relaxation is as difficult to solve to optimality as the MINLP}.

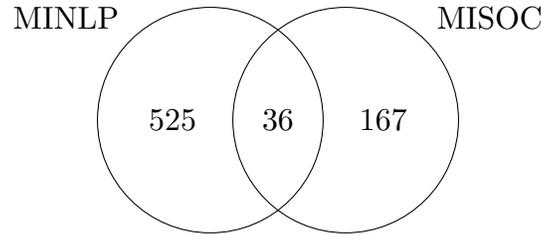
\begin{figure}
\centering
\begin{tikzpicture}
\node [draw,
    circle,
    minimum size =3cm,
    label={135:MINLP}] (A) at (0,0){};
\node at (-0.5, 0) {525};
\node [draw,
    circle,
    minimum size =3cm,
    label={45:MISOC}] (B) at (1.8,0){};
\node at (2.3, 0) {167};
\node at (0.9,0) {36};
\end{tikzpicture}
\caption{Venn diagram of timed-out instances for GasLib-4197 network.}
    \label{fig:4197-venn}
\end{figure}

The Venn diagram in Fig. \ref{fig:4197-venn} shows the number of instances that timed-out for the MINLP and the MISOC formulation. Observe that there are $36$ instances for which the MINLP was solved to global optimality but the MISOC relaxation terminated without computing a feasible solution at the end of $1000$ seconds. We also note that, the linear relaxation for all these $36$ instances converged within a few seconds and closed the gap on all of them (see Table \ref{tab:gaslib-4197-stats}, third row). The Table \ref{tab:gaslib-4197-stats} also shows the statistics of relative gap for all the instances for which the MINLP was successfully solved to global optimality. 

\begin{table}[h!]
    \centering
    \caption{Statistics of relative gaps in (\%) for GasLib-4197 instances for which the MINLP and one of the relaxations converged to its optimal solution within the computation time limit.}
    \label{tab:gaslib-4197-stats}
    \begin{threeparttable}
    \begin{tabular}{ccccc}
         \toprule 
         formulation & minimum & maximum & mean & std. dev. \\ 
         \midrule 
         Linear relaxation\tnote{a} & 0.00 & 0.54 & 0.00 & 0.01 \\ 
         MISOC relaxation\tnote{a} & 0.00 & 0.55 & 0.00 & 0.01 \\ 
         Linear relaxation\tnote{b} & 0.00 & 0.00 & 0.00 & 0.00 \\ 
         \bottomrule 
    \end{tabular}
    \begin{tablenotes}\footnotesize
        \item[a] instances for which the MISOC converged to its optimal solution
        \item[b]instances for which the MISOC timed out
    \end{tablenotes}
    \end{threeparttable}
\end{table}

\subsection{Summary of the results} \label{subsec:summary}
In summary, out of the $1234$ instances for GasLib-134, both the relaxations were able to prove optimality in $1232$ instances and were able to provide a certificate of infeasibility for the remaining $2$ instances. For the $4227$ instances of the GasLib-582 network, the MISOC relaxation closed the gap on $4212$ instances and obtained a relaxed solution that was within $1\%$ of the globally optimal solution to the OGF for $4215$ instances, whereas for the linear relaxation, these numbers were $4196$ and $4215$, respectively. Both the relaxations were ineffective in detecting infeasibilities in the $7$ instances for which the MINLP was infeasible. And finally, for the $2859$ instances of GasLib-4197 network, the MISOC relaxation closed the gap on $2025$ instances and obtained a relaxed solution that was within $1\%$ of the globally optimal solution for $2058$ instances; these numbers for the linear relaxation were $2181$ and $2225$, respectively. Similar to the instances corresponding to GasLib-582, both relaxations failed to detect infeasibilites in the $68$ infeasible instances. In total, the relaxations were able to compute the global optimal solution or a solution within $1\%$ of the global optimal solution to OGF for $7672/8320$ instances, they provided a certificate of infeasibility for $2/8320$ instances, were not able to provide a certificate of infeasibility when the MINLP was infeasible for $75/8320$ instances. For the remaining instances, the MINLP either timed out without providing a feasible solution or errored out. 
 
\section{Concluding Remarks} \label{sec:conclusion}
In this paper, we presented a novel MINLP formulation for the OGF problem for a non-ideal gas. A polyhedral relaxation is proposed for the relaxing the nonlinear terms that arise in the mathematical model of the physical components due to the non-ideal equation of state. This polyhedral relaxation is in turn used to construct a novel linear and an MISOC relaxation for the MINLP with the latter being based on existing work in the literature. Extensive computational results on benchmark gas network instances show the efficacy of the proposed relaxations. In particular, the linear relaxation produced a solution within $1\%$ of the global optimal solution with a few seconds of computation time even on very large networks at the scale that was not done before in the literature. Future work would focus on extending these relaxations for the dynamic version of the OGF, that model the transient behavior of gas flow. 

\ACKNOWLEDGMENT{%
The authors acknowledge the funding provided by LANL’s Directed Research and Development (LDRD) project: ``20220006ER: Fast, Linear Programming-Based Algorithms with Solution Quality Guarantees for Nonlinear Optimal Control Problems''. The research work conducted at Los Alamos National Laboratory is done under the auspices of the National Nuclear Security Administration of the U.S. Department of Energy under Contract No. 89233218CNA000001.
}

%
%
%


\bibliographystyle{informs2014} 
\bibliography{references} 


\end{document}